\documentclass[twoside,11pt,leqno]{article}
\usepackage{amssymb}
\usepackage{amsthm}
\usepackage{srcltx}
\usepackage[tbtags]{amsmath}
\usepackage{doc}
\usepackage{latexsym}
\usepackage{eucal}
\usepackage{enumerate}
\usepackage{amssymb}
\usepackage{mathrsfs}
\usepackage{amsmath,amssymb,amsthm,latexsym,amscd,mathrsfs}
\usepackage{indentfirst}
\usepackage{stmaryrd}
\usepackage{graphicx}

\theoremstyle{plain}
\newtheorem*{thm A}{Theorem~A}
\newtheorem*{thm B}{Theorem~C}
\newtheorem*{thm C}{Theorem~D}
\newtheorem*{Main Theorem}{Main Theorem}
\newtheorem*{pro A}{Proposition~E}
\newtheorem*{pro B}{Proposition~F}
\newtheorem*{proposition A}{Theorem~B}

\numberwithin{equation}{section}
\newtheorem{theorem}{Theorem}[section]

\newtheorem{lemma}[theorem]{Lemma}
\newtheorem{proposition}[theorem]{Proposition}
\newtheorem{remark}[theorem]{Remark}

\theoremstyle{definition}

\newcommand{\psum}{{+_{\negthinspace\kern-2pt p}}\,}

\begin{document}
\small{\addtocounter{page}{0} \pagestyle{plain}


\renewcommand{\baselinestretch}{1.5}
\renewcommand{\theequation}{\thesection.\arabic{equation}}
 \makeatletter
 \renewcommand{\section}{\@startsection {section}{0}{0mm}
                        {\baselineskip}{0.3\baselineskip}{\normalfont\bfseries}}
                        \makeatother
                       \makeatother
\makeatletter
   \renewcommand{\subsection}{\@startsection {subsection}{0}{0mm}
                        {\baselineskip}{0.3\baselineskip}{\normalfont\bfseries}}
                        \makeatother

\noindent{\large\bf  Isoparametric theory and its applications}
\vspace{0.15in}\\
\noindent{\sc Zizhou Tang}
\newline
{\small Chern Institute of Mathematics, Nankai University, Tianjin 300071, P. R. China;\\
{\small e-mail} : { \verb|zztang@nankai.edu.cn| }
\vspace{3mm}

\noindent{\sc Wenjiao Yan}
\newline
{\small School of Mathematical Sciences, Beijing Normal University, Beijing 100875, P. R. China;\\
{\small e-mail} : { \verb|wjyan@bnu.edu.cn| }
\vspace{3mm}
%



{\footnotesize
\noindent Abstract.
This is a survey on the recent progress in several applications of isoparametric theory, including an affirmative answer to Yau's conjecture on the
first eigenvalue of Laplacian in the isoparametric case, a negative answer to Yau's 76th problem in his Problem Section, new examples of Willmore submanifolds in spheres,
a series of examples to Besse's problem on the generalization of Einstein condition, isoparametric functions on exotic spheres, counterexamples to two
conjectures of Leung, as well as surgery theory on isoparametric foliation.}


\section {Preliminaries}
\setcounter{equation}{0}
\renewcommand{\theequation}{1.\arabic{equation}}
\vspace{2mm}

The modern isoparametric theory originates from E. Cartan. After the fundamental work of M\"{u}nzner, the isoparametric theory has become more and more fascinating to geometers.
Research on classifications and applications of isoparametric foliation (isoparametric functions, isoparametric hypersurfaces) in spheres have been very active in recent decades (cf. \cite{CR2}).
For classifications, see \cite{CCJ}, \cite{Imm}, \cite{Chi1, Chi2}, \cite{DN}, and \cite{Miy,Miy2}; for applications, see for example, \cite{GT1, GT2, GT3}, \cite{GTY}, \cite{GQ}, \cite{LY},
\cite{QT1, QT2, QT3}, \cite{TXY}, \cite{TY1, TY2, TY3, TY4}, \cite{Xie}. In a Riemannian manifold, an isoparametric hypersurface is a hypersurface whose nearby parallel hypersurfaces all have constant
mean curvatures. While in real space forms, an equivalent definition of isoparametric hypersurface given by E. Cartan is that the hypersurface has constant principal curvatures.

In this paper, we are mainly concerned with isoparametric foliation in the unit sphere $S^{n+1}(1)$. We start by giving a definition of isoparametric functions. A smooth function $f$ on $S^{n+1}(1)$ is called \emph{isoparametric}, if it satisfies
\begin{eqnarray}{}
 |\nabla f|^2 = ~b(f),\label{isop1}\\
 ~\Delta f ~=~a(f), \label{isop2}
 \end{eqnarray}
where $b$ and $a$ are smooth and continuous functions on $\mathbb{R}$, respectively. A function satisfies only (\ref{isop1}) is called \emph{transnormal}. The geometric meaning of (\ref{isop1}) and (\ref{isop2}) is that
the regular level hypersurfaces of $f$ are parallel with each other and have constant mean curvatures. In this sense, the regular level hypersurfaces of $f$
are called \emph{isoparametric hypersurfaces}, and the two singular level sets of $f$ are called \emph{focal submanifolds}, which are denoted by $M_+$, $M_-$ according to
they are inverse images of maximum or minimum values of $f$, respectively. Nomizu \cite{Nom} showed that $M_{\pm}$ are minimal submanifolds in $S^{n+1}(1)$.
For the isoparametric theory on a general Riemannian manifold, see \cite{Wan}, \cite{GT1}.

Let $M^n$ be a closed isoparametric hypersurface in $S^{n+1}(1)$. Denote by $g$ the number of its distinct principal curvatures $k_i$ ($k_1>\cdots >k_g$) with multiplicity $m_i$($i=1,...,g$). A fundamental result of M{\"u}nzner \cite{Mun} states that $g\in\{1,2,3,4,6\}$ and $m_i=m_{i+2}$ (subscripts mod $g$), furthermore, $M^n$ can be realized as a level hypersurface in $S^{n+1}$ of
a homogeneous polynomial $F$ of degree $g$ on $\mathbb{R}^{n+2}$
satisfying \emph{Cartan-M\"{u}nzner equations}:
\begin{eqnarray}
|\nabla F|^2&=&g^2|x|^{2g-2}, \nonumber\\
 \Delta F&=&\frac{g^2}{2}(m_2-m_1)|x|^{g-2}. \nonumber
\end{eqnarray}
The function $f=F~|~_{S^{n+1}(1)}$ takes values in $[-1, 1]$. As defined before, $M_+=f^{-1}(1)$ and $M_-=f^{-1}(-1)$. They are minimal submanifolds in $S^{n+1}(1)$, with codimensions
$m_1+1$ and $m_2+1$, respectively.

When $g \leq3$, the classification for isoparametric hypersurfaces are accomplished by E. Cartan, who proved them to be homogeneous, namely, they can be characterized as the principal orbits of the isotropy representation of some rank $2$ symmetric spaces;
when $g=4$, after Abresch \cite{Abr}, Tang \cite{Tan}, Fang \cite{Fan} and Stolz \cite{Sto} determined all the possible multiplicity pairs $(m_1, m_2)$, Cecil-Chi-Jensen \cite{CCJ}, Immervoll \cite{Imm} and Chi \cite{Chi1, Chi2} showed that the isoparametric hypersurfaces are either of OT-FKM type, or homogeneous with $(m_1, m_2)=(2, 2), (4, 5)$;
when $g=6$, Abresch \cite{Abr} showed that the multiplicities $m_1=m_2= 1$ or $2$, and then Dorfmeister-Neher \cite{DN} and Miyaoka \cite{Miy, Miy2} showed that they are homogeneous.

Therefore, all the inhomogeneous isoparametric hypersurfaces in $S^{n+1}(1)$ appear in the OT-FKM type. Now we turn to define it. For a symmetric Clifford system $\{P_0,\cdots,P_m\}$ on
$\mathbb{R}^{2l}$, \emph{i.e.},
\begin{equation*}
 P_i=P^t_i , ~~P_iP_j+P_jP_i=2\delta_{ij}I_{2l},
\end{equation*}
following Ozeki and Takeuchi \cite{OT}, Ferus, Karcher and M\"{u}nzner \cite{FKM} constructed a homogeneous polynomial $F$
of degree $4$ on $\mathbb{R}^{2l}$:
  \begin{equation*}
        F(x) = |x|^4 - 2\sum_{i = 0}^{m}{\langle P_ix,x\rangle^2},
  \end{equation*}
where $l=k\delta(m)$, $k$ is a positive integer, $\delta(m)$ is the dimension of an irreducible module of the Clifford algebra $\mathcal{C}_{m-1}$ which is valued:
\begin{center}
\begin{tabular}{|c|c|c|c|c|c|c|c|c|c|}
\hline
$m$ & 1 & 2 & 3 & 4 & 5 & 6 & 7 & 8 & $\cdots$ $m$+8 \\
\hline
$\delta(m)$ & 1 & 2 & 4 & 4 & 8 & 8 & 8 & 8 & ~16$\delta(m)$\\
\hline
\end{tabular}
\end{center}
\cite{FKM} showed that the polynomial $F$ satisfies the Cartan-M\"{u}nzner equations with multiplicity pairs $(m_1, m_2)=(m, l-m-1)$.

\section{Yau's conjecture on the first eigenvalue}

Let $M^n$ be a closed Riemannian manifold. It is well known that the Laplace operator $\Delta$ is an elliptic operator and has a discrete spectrum:
$$\{0=\lambda_0(M)<\lambda_1(M)\leqslant \lambda_2(M)\leqslant \cdots \leqslant\lambda_k(M),\cdots, \uparrow \infty\}.$$
When $M^n$ is an $n$-dimensional minimal immersed submanifold in a unit sphere $S^N$, Takahashi \cite{Tak} proved that the coordinate functions are eigenfunctions of $\Delta$ with
corresponding eigenvalue $n$.
Thus the multiplicity of $n$ is at least $N+1$ and $\lambda_1(M^n)\leq n$.
In this regard, S.T.Yau raised the following problem in his problem section in 1982:
\vspace{2mm}

\noindent
\emph{\textbf{Conjecture(\cite{Yau})}: Let $M^n$ be a closed embedded minimal hypersurface in the unit sphere $S^{n+1}(1)$. Then
$$\lambda_1(M^n)=n.$$}

Certainly, this is a conjecture of cardinal significance. For example, a well-known conjecture of Lawson, which states that
the only embedded minimal torus in $S^3(1)$ is the Clifford torus, was proved by Brendle \cite{Bre}. While in
1986, Montiel and Ros \cite{MR} showed that
Lawson's conjecture is just the torus case of Yau's conjecture on the first eigenvalue.

Attacking Yau's conjecture, Choi and Wang \cite{CW} proved that $\lambda_1(M^n)\leq n/2$.
We consider a restricted version of this problem---the isoparametric hypersurface, which is naturally an embedded hypersurface.
Muto, Ohnita, Urakawa \cite{MOU} and Kotani \cite{Kot} considered some homogeneous isoparametric hypersurfaces case.
Combining their results with the classification, one concludes that the first eigenvalues of minimal isoparametric hypersurfaces with $g=1,2,3,6$
in $S^{n+1}(1)$ equal to the dimension $n$. Furthermore, Muto \cite{Mut} also obtained the same conclusion for a small number of inhomogeneous isoparametric hypersurfaces with $g=4$.
Independent of homogeneity, we proved the following wide generalized result:
\begin{theorem}[\cite{TY2}]
Let $M^n$ be a minimal isoparametric hypersurface in $S^{n+1}(1)$ with $g=4$ and $m_1, m_2\geqslant 2$. Then
$$\lambda_1(M^n)=n$$
with multiplicity $n+2$.
\end{theorem}

Takagi \cite{Tak1} asserted that the isoparametric hypersurface with $g=4$ and multiplicities $(1, k)$ must be homogeneous.
As mentioned before, this case was proved by \cite{MOU}. Therefore, we finally arrive at an affirmative answer to Yau's conjecture on the first eigenvalue in the isoparametric case:
\begin{theorem}[\cite{TY2}]
Let $M^n$ be a minimal isoparametric hypersurface in $S^{n+1}(1)$. Then
$$\lambda_1(M^n)=n.$$
\end{theorem}

\begin{remark}
S. S. Chern conjectured that a closed, minimally immersed hypersurface in $S^{n+1}(1)$, whose second fundamental form has constant length, is isoparametric (cf.
\cite{GT3}). If this conjecture is proven, we would have settled Yau's conjecture for the
minimal hypersurfaces whose second fundamental form have constant length, which gives
us more confidence in Yau's conjecture.
\end{remark}

Since the focal submanifolds are minimal submanifolds in spheres, it would be natural to ask that
\emph{what are the first eigenvalues of the focal submanifolds $M_+$ and $M_-$?}
For the cases $g=1,2,3$, it can be trivially observed that the first eigenvalues are equal to their dimensions. For $g=4$, we proved
\begin{theorem}[\cite{TY2}]
Let $M_+$ be a focal submanifold of codimension $m_1+1$ in $S^{n+1}(1)$ with $g=4$.
If $\dim M_+\geqslant \frac{2}{3}n+1$, then
$$\lambda_1(M_+)=\dim M_+=m_1+2m_2$$
with multiplicity $n+2$. A similar conclusion holds for $M_-$ with $m_1$, $m_2$ exchanging.
\end{theorem}

In a following paper, Tang, Xie, and Yan \cite{TXY2} also proved Yau's conjecture for the minimal isoparametric hypersurface in $S^{n+1}(1)$ with $g=6$ and $m_1=m_2=2$ without using the classification,
and showed that one focal submanifold in this case has dimension as the first eigenvalue with multiplicity $n+2$.

Based on all the results above, we would like to ask the following question:
\vspace{2mm}

\noindent
\emph{\textbf{Problem (\cite{TY2})}. Let $M^d$ be a closed minimal embedding submanifold in $S^{n+1}(1)$. If the
dimension $d$ of $M^d$ satisfies $d\geqslant \frac{2}{3}n+1$, then is it true that $\lambda_1(M^d)=d$?
}
\vspace{2mm}

To the best of our knowledge, there are no counterexamples and the dimensional condition is very necessary.
For instance, Solomon \cite{Sol} found some eigenfunctions on the focal submanifold $M_-$ of OT-FKM type,
which correspond to eigenvalues $4m_1$. However, if we require $m_1<\frac{1}{2}m_2$ ( for example
$(m_1, m_2)=(2, 5)$), which implies that $d< \frac{2}{3}n+1$, then $\lambda_1(M_-)\leqslant 4m_1<\dim M_-=2m_1+m_2.$
This shows that the assumption on the dimension of the submanifold in Theorem 2.4 is essential.

\section{Negative answer to Yau's 76th problem}
There are few manifolds whose eigenvalues are clearly known, not to mention the eigenfunctions.
The numbers of critical points of eigenfunctions are even more difficult to determine. However, as
S.T.Yau pointed out, this number is closely related to many important questions, which
makes it worthy of being studied extensively. In this regard, S.T.Yau \cite{Yau} raised the following
problem as the 76th problem in his problem section:
\vspace{2mm}

\noindent
\emph{\textbf{Conjecture(\cite{Yau})}. The number of critical points of the $k$-th eigenfunction on a
compact Riemannian manifold increases with $k$.}
\vspace{2mm}

Yau also investigated this problem in the surface case (cf. \cite{Yau2}).
In 1999, Jakobson and Nadirashvili \cite{JN} constructed a metric on a $2$-dimensional
torus and a sequence of eigenfunctions such that the corresponding eigenvalues go to
infinity while the number of critical points remains $16$. However, their example does not deny Yau's conjecture in the
sense of ``non-decreasing".

Taking advantage of isoparametric theory, we constructed other examples to rule out the possibility of ``non-decreasing".

\begin{theorem}[\cite{TY3}]
Let $M^n$ be the minimal isoparametric hypersurface of OT-FKM type in $S^{n+1}(1)$. Then there exist three eigenfunctions $\varphi_1$, $\varphi_2$ and $\varphi_3$ defined
on $M^n$, corresponding to eigenvalues $n$, $2n$ and $3n$, whose critical sets consist of $8$
points, a submanifold and $8$ points, respectively. For specific, $\varphi_1$ and $\varphi_3$ are both Morse
functions; $\varphi_2$ is an isoparametric function on $M^n$, whose critical set $C(\varphi_2)$ is:
$C(\varphi_2)=N_+\cup N_-$, $\dim N_+ = \dim N_-= n-m ~(1 \leq m<n)$, where the number $m$ the same with that in the definition of OT-FKM type.
\end{theorem}

\begin{remark}
The Morse number (the minimal number of critical points of all Morse
functions) of a compact isoparametric hypersurface with $g=4$ in $S^{n+1}(1)$ is equal to $2g=8$ (cf. \cite{CR}).
\end{remark}

For specific, choosing a point $q_1\in S^{n+1}(1)\backslash \{M_+, M_-, M^n\}$, we defined $\varphi_1$, $\varphi_2$ and
$\varphi_3$ as follows:
\begin{eqnarray*}
&& \varphi_1: M^n\rightarrow \mathbb{R},\qquad\qquad \varphi_2: M^n\rightarrow \mathbb{R}\qquad\qquad \varphi_3: M^n\rightarrow \mathbb{R}\nonumber\\
&& \qquad x \mapsto \langle x, q_1\rangle, \qquad\qquad\quad x \mapsto \langle Px, x\rangle \quad\qquad\qquad x \mapsto \langle \xi(x), q_1\rangle
\end{eqnarray*}
where $\xi$ is a unit normal vector field on $M^n$,
and $P\in \Sigma (P_0,...,P_m)$, the unit sphere in
$\mathrm{Span}\{P_0,...,P_m\}$, which is called \emph{the Clifford sphere}.

Moreover, denoting by $c_0$ the regular value corresponding to $M^n$, we showed that $\varphi_2$ satisfies
\begin{equation*}\label{isoparametric varphi2}
\left\{ \begin{array}{ll}
|\nabla \varphi_2|^2=4(1-\frac{2}{1-c_0}\varphi_2^2)\\
\quad\Delta \varphi_2 = -2n \varphi_2.
\end{array}\right.
\end{equation*}
Furthermore, we proved that the focal submanifolds $N_{\pm}$ of $\varphi_2$ on $M^n$ are both diffeomorphic to $M_+$ of OT-FKM type.

In a similar way, we also constructed another counterexample on the focal submanifold $M_-$ of OT-FKM type.

\section{Willmore, Einstein properties of focal submanifolds and Besse's problem on the generalization of Einstein condition}
Let $x:M^n \rightarrow S^{n+p}(1)$ be an immersion from an
$n$-dimensional compact manifold to an $(n + p)$-dimensional unit
sphere. Denote by $S$ the norm square of the second fundamental form, and $H$ the norm of the mean curvature
vector. Then $M^n$ is called a \emph{Willmore submanifold} in $S^{n+p}(1)$ if it is an extremal submanifold of the
Willmore functional, which is a conformal invariant (\emph{cf.}
\cite{PW, Wan'}):
$$W(x)=\int_{M^n} (S-nH^2)^{\frac{n}{2}}dv.$$

In conjunction with the fact that the focal submanifolds of an
isoparametric hypersurface in $S^{n+1}(1)$ are minimal submanifolds of
$S^{n+1}(1)$ with constant $S$, we established the following

\begin{theorem}[\cite{TY1, QTY}]
{\itshape Both focal submanifolds of every isoparametric
hypersurface in unit spheres with $g=4$ are Willmore submanifolds in $S^{n+1}(1)$.}
\end{theorem}

Moreover, as a generalization, Xie reproved our result and showed

\begin{theorem}[\cite{Xie}]
{\itshape Both focal submanifolds of every isoparametric hypersurface in
unit spheres with $g=6$ are Willmore in $S^{n+1}(1)$.}
\end{theorem}

It is well known that an Einstein manifold minimally immersed in a unit sphere is Willmore. So we
are interested in whether focal submanifolds are Einstein. After a deeper exploration on the shape operator of OT-FKM type and the homogeneous cases, we obtained

\begin{theorem}[\cite{QTY}]
{\itshape For the focal submanifolds of an isoparametric
hypersurface in $S^{n+1}(1)$ with four distinct principal curvatures,
we have
\begin{itemize}
\item[(i)] All the $M_-$ of OT-FKM type are not Einstein; the $M_+$ of
           OT-FKM type is Einstein if and only if it is diffeomorphic to $Sp(2)$ in the homogeneous case with
           $(m_1,m_2)=(4,3)$.
\item[(ii)] In the case $(m_1,m_2)=(2,2)$, the focal submanifold
           diffeomorphic to $\widetilde{G}_2(\mathbb{R}^5)$ is Einstein, while the other one
           diffeomorphic to $\mathbb{C}P^3$ is not.
\item[(iii)]In the case $(m_1,m_2)=(4,5)$, both focal submanifolds are not Einstein.
\end{itemize}
}
\end{theorem}

\begin{remark}
In the OT-FKM type, there are two incongruent examples corresponding to $(m_1,m_2)=(4,3)$, one is homogeneous
while the other is not. It is surprising that only the focal submanifold $M_+$ of the homogenous case is Einstein.
\end{remark}

Noticing that there are only two Einstein manifolds among all the focal submanifolds with $g=4$, we are curious about their Einstein-like properties.
Let $\mathcal{E}$ denote the set of Einstein manifolds, $\mathcal{P}$ the set of Riemannian manifolds with parallel Ricci tensor, and $\mathcal{S}$
the set of Riemannian manifolds with constant scalar curvatures. It is obvious that $\mathcal{E}\subset \mathcal{P}\subset \mathcal{S}$.
As further generalizations of the Einstein condition, A. Gray \cite{Gra} introduced two significant classes of Einstein-like Riemannian manifolds, that is, $\mathcal{A}$-manifolds, on which the Ricci tensor $\rho$
is cyclic parallel: $\nabla_i\rho_{jk}+\nabla_j\rho_{ki}+\nabla_k\rho_{ij}=0$, and $\mathcal{B}$-manifolds, on which the Ricci tensor $\rho$ is a Codazzi tensor: $\nabla_i\rho_{jk}-\nabla_j\rho_{ik}=0$.
These two classes have been investigated extensively since then.
Gray also showed that $\mathcal{A}$-manifolds and $\mathcal{B}$-manifolds both have constant scalar curvatures, and the intersection of them is excatly $\mathcal{P}$.
More precisely, $\mathcal{A}$ and $\mathcal{B}$ are the only classes between $\mathcal{P}$ and $\mathcal{S}$ from the
view of group representations.

The examples of $\mathcal{A}$-manifolds are not rare in the literature. For instance, the D'Atri spaces (Riemannian manifolds with
volume preserving geodesic symmetries) belong to the class $\mathcal{A}$. However, the known examples are mostly (locally) homogeneous.
In this regard, Besse (\cite{Bes}, 16.56(i), pp.451) posed the following challenging problem as one of ``some open problems" :

\vspace{2mm}

\noindent
\emph{\textbf{Problem (\cite{Bes})}. Find
examples of $\mathcal{A}$-manifolds, which are neither locally
homogeneous, nor locally isometric to Riemannian products and have
non-parallel Ricci tensor}.
\vspace{2mm}

Aiming for this problem, Jelonek \cite{Jel} and Pedersen, Tod \cite{PT}
constructed $\mathcal{A}$-manifolds on $S^1$-bundles over locally
non-homogeneous K\"{a}hler-Einstein manifolds, and on $S^1$-bundles
over a $K3$ surface, from defining Riemannian submersion metric on
the $S^1$-bundles. But their examples are not simply-connected.

On the other hand, according to \cite{QTY}, there are only two Einstein manifolds among all the focal submanifolds
of $g=4$. Thus a further and natural question arises:

\emph{Are the focal submanifolds with $g=4$ Ricci parallel, $\mathcal{A}$-manifolds, or
$\mathcal{B}$-manifolds ? }

Starting from this point, we found a series of
simply-connected examples with natural metric for this open problem of Besse.
In fact, we firstly showed

\begin{theorem}[\cite{TY4, LY}]
All the focal submanifolds of isoparametric hypersurfaces in spheres
with $g=4$ are $\mathcal{A}$-manifolds.
\end{theorem}

Then comparing with the classification of Einstein manifolds for the focal submanifolds with $g=4$, we obtained
the classification of Ricci parallel manifolds for the focal submanifolds with $g=4$ as follows:

\begin{theorem}[\cite{TY4}]
For the focal submanifolds of isoparametric hypersurfaces in spheres
with $g=4$, we have
\begin{itemize}
\item [(i)] The $M_+$ of OT-FKM type is Ricci parallel if and only if
$(m_1,m_2)=(2,1), (6,1)$, or it is diffeomorphic to $Sp(2)$ in the
homogeneous case with $(m_1,m_2)=(4,3)$; while the $M_-$ of OT-FKM
type is Ricci parallel if and only if $(m_1,m_2)=(1,k)$.

\item[(ii)]For $(m_1,m_2)=(2,2)$, the one
           diffeomorphic to $\widetilde{G}_2(\mathbb{R}^5)$ is Ricci parallel, while the other
           diffeomorphic to $\mathbb{C}P^3$ is not.
\item[(iii)]For $(m_1,m_2)=(4,5)$, both are not Ricci parallel.
\end{itemize}
\end{theorem}

Furthermore, we showed
\begin{proposition}[\cite{TY4}]
The focal submanifolds of isoparametric hypersurfaces in spheres with $g=4$
and $m_1, m_2>1$ are not Riemannian products.
\end{proposition}

Among all the focal submanifolds with $g=4$, we took one series of focal submanifolds to investigate the homogeneity:
\begin{proposition}[\cite{TY4}]
The focal submanifolds $M_+$ of OT-FKM type with $(m_1, m_2)=(3, 4k)$
are not intrinsically homogeneous.
\end{proposition}

At last, recalling the following lemma,

\begin{lemma}[\cite{Tan}]
If $m_1>1$ (resp. $m_2>1$), the focal submanifold $M_-$ (resp.
$M_+$) is simply-connected.
\end{lemma}
\noindent
we finally established
\begin{theorem}[\cite{TY4}]
The focal submanifolds $M_+$ of OT-FKM type with $(m_1, m_2)=(3, 4k)$, $k \geq 1$, are simply-connected examples to Besse's problem.
\end{theorem}

\section{Other applications and generzlizations}
\subsection{Isoparametric functions on exotic spheres}

In a general Riemannian manifold, an isoparametric function can be defined in the same way with (\ref{isop1}), (\ref{isop2}).
Moreover, Q. M. Wang\cite{Wan}, Ge-Tang\cite{GT1} 
generalized Nomizu's result about the minimality of focal submanifolds in $S^{n+1}(1)$.

\begin{theorem}[\cite{Wan},\cite{GT1}]
The focal sets $M_{\pm}$ of an isoparametric function $f$ on a complete Riemannian manifold $N$ are minimal submanifolds.
\end{theorem}

\cite{GT2} initiated the study of isoparametric functions on exotic spheres (i.e., they are homeomorphic but not diffeomorphic to the standard spheres),
and used it to attack the so-called ``smooth Poincar\'{e} conjecture in dimension four": \emph{Whether an exotic $4$-sphere exists?}
In fact, they showed

\begin{theorem}[\cite{GT1}]
Suppose $\Sigma^4$ is a homotopy $4$-sphere and it admits a transnormal
function under some metric. Then $\Sigma^4$ is diffeomorphic to the standard $S^4$.
\end{theorem}

\cite{GT1} also constructed 
explicitly a properly transnormal but not isoparametric function on the Gromoll-Meyer sphere with two points as the focal sets. It intrigued one to
study the types of focal sets in a complete Riemannian manifolds. According to \cite{Wan},
transnormal functions are necessarily Morse-Bott functions. Qian and Tang \cite{QT2}
gave the following fundamental construction of isoparametric functions

\begin{theorem}[\cite{QT2}]
Let $f$ be a Morse-Bott function on a closed manifold $N$, s.t. the critical
set has two components and each component has codimension at least $2$
in $N$. Then $N$ admits a metric, s.t. $f$ is isoparametric.
\end{theorem}

By a remarkable theorem of Smale, the above theorem has a consequence:
\begin{theorem}[\cite{QT2}]
Every homotopy $n$-sphere with $n>4$ admits a metric and an
isoparametric function with $2$ points as the focal sets.
\end{theorem}

Actually, Ge \cite{Ge} found the following correspondence:

\begin{theorem}[\cite{Ge}]
Up to foliated diffeomorphism, there exists a $1$-$1$ correspondence
between isoparametric foliations on the standard $S^n$ and on any homotopy sphere
$\Sigma^n$.
\end{theorem}

In addition, \cite{QT2} proved the following non-existence result:
\begin{theorem}[\cite{QT2}]
Every odd dimensional exotic sphere admits no totally isoparametric
functions with $2$ points as the focal set.
\end{theorem}

Here, a totally isoparametric function (defined in \cite{GTY}) is an isoparametric function so that each
regular level hypersurface has constant principal curvatures.

\subsection{New examples of harmonic maps}
Inspired by the OT-FKM construction, for a symmetric
Clifford system $\{P_0,\cdots, P_m\}$ on $\mathbb{R}^{2l}$, Qian and Tang \cite{QT3} defined
$M_i:=\{x\in S^{2l-1}(1)~|~ \langle P_0x, x\rangle = \langle P_1x, x\rangle=\cdots=\langle P_ix, x\rangle=0\}$ for $0\leq i\leq m-1$, and
obtained a minimal isoparametric hypersurface sequence
$$M_m=M_+\subset M_{m-1}\subset \cdots\subset M_0\subset S^{2l-1}(1).$$
The naturally define function $f_i: M_i\rightarrow \mathbb{R}$ by $f_i(x)=\langle P_{i+1}x,x\rangle$ for $x\in M_i$ with $Im(f_i)=[-1,1]$ is actually
an isoparametric function on $M_i$.

Similarly, by defining $N_i:=\{x\in S^{2l-1}(1)=\langle P_0x, x\rangle^2+\langle P_1x, x\rangle^2+\cdots+\langle P_ix, x\rangle^2=1\}$ for $2\leq i\leq m$, they
constructed another minimal isoparametric hypersurface sequence
$$N_1\subset N_2\subset \cdots\subset N_m=M_-\subset S^{2l-1}(1).$$
 The naturally define function $g_i: N_i\rightarrow \mathbb{R}$ by $g_i(x)=\langle P_{i}x,x\rangle$ for $x\in N_i$ with $Im(g_i)=[-1,1]$ is actually
an isoparametric function on $N_i$.

Based on these, they gave estimates on eigenvalues of the Laplacian of the focal submanifolds of isoparametric hypersurfaces in unit spheres, improving
results of \cite{TY2} and \cite{TXY2}.

Let $N$ be a closed Riemannian manifold.
A map $\varphi : N \rightarrow S^n(1)$ is called an \emph{eigenmap} if the $\mathbb{R}^{n+1}$-components are eigenfunctions of the Laplacian of $N$ and all have
the same eigenvalue. In particular, $\varphi$ is a \emph{harmonic map}, i.e., a critical point of the energy functional defined by $E(f)=\frac{1}{2}\int_M|df|^2dV_M$.  In 1980, Eells and Lemaire posed the following problem.
\vspace{2mm}

\noindent
\emph{\textbf{Problem (\cite{EL})}.
Characterize those compact $N$ for which there is an eigenmap $\varphi : N \rightarrow S^n(1)$ with $\dim(N)\geq n$?}
\vspace{2mm}

As an application of their construction, \cite{QT3} showed
\begin{theorem}[\cite{QT3}]
Let $\{P_0,\cdots,P_m\}$ be a symmetric Clifford system on $\mathbb{R}^{2l}$.

\begin{itemize}
  \item [(1)]  For $0\leq i\leq m-1$, both of the focal maps $\varphi_{\pm\pi/4}: M_{i+1}\rightarrow \mathcal{U}_{\pm 1}=f_i^{-1}(\pm 1)\cong S^{l-1}(1)$ defined by
  $\varphi_{\pm\pi/4}(x)=\frac{1}{\sqrt{2}}(x\pm P_{i+1}x)$ for $x\in M_{i+1}$ are submersive eigenmaps with the same eigenvalue $2l-i-3$.
 \item [(2)] For $2\leq i\leq m$, both of the focal maps $\phi_{\pm\pi/4}: N_{i-1}\rightarrow \mathcal{V}_{\pm 1}=g_i^{-1}(\pm 1)\cong S^{l-1}(1)$  defined by
$\phi_{\pm\pi/4}(x)=\frac{1}{\sqrt{2}}(x\pm P_{i}x)$ for $x\in N_{i}$ are submersive eigenmaps with the same eigenvalue $l+i-2$.

\end{itemize}
\end{theorem}

Meanwhile, they also considered the case for isoparametric foliations on unit spheres and obtained

\begin{proposition}[\cite{QT3}]
Let $M$ be a closed isoparametric hypersurface in a unit sphere. Then
every focal map from $M$ to its focal submanifolds $M_+$ or $M_-$ is harmonic.
\end{proposition}

\subsection{Counterexamples to Leung's two conjectures}
Let $W^n$ be a closed Riemannian manifold minimally immersed
in $S^{n+p}(1)$. Define an extrinsic quantity $\sigma(W)=\max\{|B(X,X)|^2~|~X\in TW, |X|=1\}$, where $B$ is the second fundamental form.
In 1986, Gauchman \cite{Gau} established a well-known rigidity theorem which states that if
$\sigma(W)<1/3$, then the submanifold $W$ must be totally geodesic. When the dimension $n$ of
$W$ is even, the rigidity theorem above is optimal. When the dimension $n$ is odd and $p > 1$, the conclusion still holds under
a weaker assumption $\sigma(W)\leq 1/(3-2/n)$. In 1991, Leung \cite{Leu} proved that if $n$ is odd, then a closed minimally immersed submanifold
$W^n$ with $\sigma(W)\leq n/(n-1)$ is totally geodesic provided that the normal connection is flat.
Based on this fact, he proposed the following conjecture.
\vspace{2mm}

\noindent
\emph{\textbf{Conjecture 1 (\cite{Leu})}.
If $n$ is odd, $W^n$ is minimally immersed in $S^{n+p}(1)$ with $\sigma(W)\leq n/(n-1)$, then $W$ is homeomorphic to $S^n$.}
\vspace{1mm}

By investigating the second fundamental form of the Clifford minimal hypersurfaces in unit
spheres, Leung also posed the following stronger:
\vspace{2mm}

\noindent
\emph{\textbf{Conjecture 2 (\cite{Leu})}.
If $n$ is odd, $W^n$ is minimally immersed in $S^{n+p}(1)$ with $\sigma(W)< (n+1)/(n-1)$, then $W$ is homeomorphic to $S^n$.}
\vspace{1mm}

For minimal submanifolds in unit spheres with flat normal connections, Conjecture 2
was proved by Hasanis and Vlachos \cite{HV}. Using isoparametric foliation, \cite{QT3}
found that

\begin{theorem}[\cite{QT3}]
The focal submanifold $M_+$ (resp. $M_-$) with $g=4$ and dimension $m_1+2m_2$ (resp. $2m_1+m_2$) is minimal in $S^{n+1}(1)$ with
$\sigma(M_{\pm})=1$. However, $M_{\pm}$ are not homeomorphic to the spheres.
\end{theorem}

Notice that if $m_1$ (resp. $m_2$) is odd, then $M_+$ (resp. $M_-$) is a counterexample
to both conjectures of Leung.

\subsection{Schoen-Yau-Gromov-Lawson's surgery theory and isoparametric foliation}
``Which manifolds admit Riemannian metrics of positive scalar curvature" has been a core issue in the study of geometry.
Let $M$ be a compact spin manifold, if it carries a Riemannian metric of positive scalar curvature, Lichnerowicz \cite{Lic} showed that
its $\widehat{A}$-genus (for $4k$ dimension) vanishes, and Hitchin \cite{Hit} generalized the result to that the invariant $\alpha(M)$ (for general dimension) vanishes, which leads to that
half of the exotic spheres in dimensions $8k+1$ and $8k+2$ cannot carry metrics of positive scalar curvature.
Another remarkable step was made by Schoen and Yau \cite{SY}, and
Gromov and Lawson \cite{GL} independently, who established the following ``surgery
theorem":
\begin{theorem}[\cite{SY, GL}]
Let $X$ be a compact manifold which carries a Riemannian metric of positive scalar curvature. Then any manifold which can be obtained
from $X$ by performing surgeries in codimension $\geq 3$ also carries a metric of positive scalar curvature.
\end{theorem}

We performed their surgery theory at the compact minimal isoparametric hypersurface $Y^{n}$ of the unit sphere $S^{n+1}(1)$, $n\geq 2$.
Here $Y^{n}$ separates $S^{n+1}$ into two ball bundles $S^{n+1}_+$ and $S^{n+1}_-$ over each of the focal submanifolds.
Gluing $S^{n+1}_+$ (resp. $S^{n+1}_-$) along $Y^n$ with itself, we obtained the double manifolds $D(S^{n+1}_+)$ (resp. $D(S^{n+1}_-)$).
For the geometrical properties of $D(S^{n+1}_{\pm})$, we proved

\begin{theorem}[\cite{TXY}]
Each of doubles $D(S^n_{\pm})$ admits a metric of positive scalar curvature. Moreover,
there is still an isoparametric foliation in $D(S^n_+)$ or $D(S^n_-)$.
\end{theorem}

Furthermore, we proved that the double manifold $D(S^{n+1}_{\pm})$ are stably parallelizable manifolds,
and studied the topology of $D(S^{n+1}_{\pm})$ for the homogeneous cases and OT-FKM type separately. In the homogeneous case, we determined completely the isotropy subgroups of the singular orbits, while in the
OT-FKM type, using K-theory, and the investigation of the associated Clifford modules, we got some classification of the homotopy types.

In \cite{PQ}, they studied Ricci curvature properties of double manifolds together with isoparametric foliations, and improved
Theorem 5.11 to positive Ricci curvature.

\end{document}